%% file: undirected2.tex
\documentclass[a4paper,11pt]{article}      

\usepackage{graphicx,color}
\usepackage{algorithmic}
\usepackage{algorithm}
\usepackage{epsfig}
\usepackage{times}

\usepackage{amsmath,amssymb}
\newcommand{\Rnn}{\R^{n \times n}}
\newcommand{\Rnm}{\R^{n \times m}}
\usepackage{amssymb,eucal,amsmath}
\usepackage{graphics,graphicx} 

\newcommand{\bullettop}[1]{\overset{\bullet}{#1}}

\newcommand{\rank}{\mbox{rank }}

\newcommand{\Cinf}{\mathfrak{C}^{\infty}}

\newcommand{\CinfRRv}{\Cinf(\R,\Rv)}

\newcommand{\Rv}{\R^{v}}

\newcommand{\Cn}{\C^{n}}
\newcommand{\Rn}{\R^{n}}

\newcommand{\Rs}{\R[s]}

\newcommand{\zeros}{\mathrm{zeros}}

\newtheorem{theorem1}{Theorem}
\newtheorem{definition}[theorem1]{Definition}
\newtheorem{corollary}[theorem1]{Corollary}
\newtheorem{prp}[theorem1]{Proposition}

\newcommand{\Rpvs}{{\R^{p\times v}[s]}}

\newcommand{\Rvvs}{{\R^{v\times v}[s]}}
\newcommand{\Rpps}{{\R^{p\times p}[s]}}

\newcommand{\Monw} {M(\der)w}
\newcommand{\deri}[1]{\mbox{$\frac{{\mathrm d}^{#1}}{{\mathrm d}t^{#1}}$}}
\newcommand{\der}{\deri{}}

\newtheorem{algorithm1}[theorem1]{Algorithm}
\newtheorem{lemma}[theorem1]{Lemma}

\newenvironment{proof1}{\noindent
          {\bf Proof.}\ }{\hfill $\Box$ \vspace*{0.1mm}\newline}
\newenvironment{proof2}{\noindent
          {\bf Proof of}}{\hfill $\Box$ \vspace*{0.1mm}\newline}
\newtheorem{exa}[theorem1]{Example}

\newcommand{\R}{{\mathbb{R}}}

\newcommand{\C}{{\mathbb{C}}}

\newcommand{\B}{{\mathfrak{B}}}

\newenvironment{example}{\begin{exa}}{\hfill $\Box$ \end{exa}}

\newcommand{\comment}[1]{#1}

\usepackage{color}

\title{Structural controllability: an undirected graph approach}
\author{Madhu N. Belur \ and 
Sivaramakrishnan Sivasubramanian\thanks{The authors are respectively 
in the Departments of Electrical Engineering and Mathematics, Indian
Institute of
Technology Bombay, Powai Mumbai 400076, India.
Corresponding author's email:
{\tt belur@ee.iitb.ac.in} and fax number: +91.22.2572.3707~.}} %
\date{}

\begin{document}

\maketitle

\begin{abstract}
This paper addresses questions regarding controllability for
`generic parameter' dynamical systems, i.e. the question whether
a dynamical system is `structurally controllable'. Unlike conventional
methods that deal with structural controllability, our approach
uses an undirected graph: the behavioral approach to modelling
dynamical systems allows this. Given a system of linear, constant
coefficient, ordinary differential equations of any order,
we formulate necessary and sufficient conditions for
controllability in terms of weights of the edges
in a suitable bipartite graph constructed from
the differential-algebraic system.  
A key notion that helps formulate the conditions is that of a `redundant
edge'. Removal of all redundant edges makes the inferring of structural
controllability a straightforward exercise.
We use standard graph algorithms as ingredients to check these conditions and 
hence obtain an algorithm to check for structural controllability.   
We provide an analysis of the running time of our algorithm.  When our 
results are applied to the familiar
state space description of a system, we obtain a novel necessary
and sufficient condition to check structural controllability for
this description.

\end{abstract}

\section{Introduction and related work} \label{sec1}

When dealing with very large dynamical systems, numerical computation is 
often not feasible.
Under the assumption of genericity of parameters, one can answer questions about
controllability and ability to achieve arbitrary pole placement
using graph theoretic tools.
These issues are typically dealt as `structural' issues
in control, see \cite{Lin:1974,Mur:1987} and the survey
paper \cite{DioComWou:03}.
While existing techniques to address structural aspects of control start
from a (possibly singular) state space representation of the system, the results in this paper
apply to more general models of dynamical systems: linear differential-algebraic equations
of possibly higher order. The behavioral theory of systems allows this general approach.  While this problem has been studied and analyzed
thoroughly since the classical paper by C.-T. Lin \cite{Lin:1974}, this
paper handles this problem using an {\em undirected graph}. This is
possible because in the behavioral approach to systems, variables
are not classified as inputs and outputs, and hence the relation between
two variables does {\em not} have to be a direction of influence
of one on another. Further, dealing with higher order differential equations
is just as easy as first order. The proposed method is very straightforward
and intuitive: we construct a weighted bipartite graph with one vertex set 
as the equations and another vertex set as the variables. 
Lack of structural controllability is shown to be equivalent to existence
of connected components of a closely related bipartite graph with some
edge-weight conditions. The non-existence of such components
can be checked using the algorithm we propose (in Section \ref{sec5}), whose
running time is quantified using standard graph algorithms.

While structural controllability is the main focus of this paper,
the results in this paper are relevant to some other questions about
generic properties of polynomial matrices.
(See Definition \ref{defn:genericity} for the precise meaning of
genericity of a property.) The first
question is under what conditions on a polynomial matrix $M$ can
we say that the {\em invariant polynomials} of $M$ are generically one.
This is nothing but the question as to when are the determinants
of all the maximal minors of a (nonsquare) polynomial matrix generically
coprime.  Another way to formulate  this question is when does the
Smith normal form of a univariate polynomial matrix generically have only ones
(and zeros) along the diagonal. Finally, this issue is equivalent
to the ability to embed a polynomial matrix as a sub-matrix of 
a unimodular matrix. (A unimodular matrix is defined in Section \ref{sec2}
as a square polynomial matrix whose determinant is a nonzero constant.)

The paper is organized as follows.  Some definitions regarding the
behavioral approach
and some graph theoretic preliminaries are covered in the following section.
Section \ref{sec3} contains some results
for square polynomial matrices. A notion called `redundant edge' is
introduced here. This notion plays a key role in this paper.
Section \ref{sec4} contains the main results of this paper: two
equivalent conditions for checking structural controllability of a dynamical
system.  Section \ref{sec5} contains an analysis into the efficiency
of the algorithm we propose for checking structural controllability.
Section \ref{sec6} specialises our results to state space systems.
In Section \ref{sec7}, we study port-terminal interconnection based models.
We also consider the standard interconnection procedures: series, parallel
and feedback and show properties of these interconnections.  In Section
\ref{sec8}, we study the situation when the polynomial matrix does not 
have full rank.  Section \ref{sec9} has some conclusive remarks.

\section{Preliminaries} \label{sec2}
The following subsection deals with polynomial matrices, the next covers
the results about the behavioral approach to modelling and control of
dynamical systems, while Subsection \ref{sub3:sec2} deals with graph
theoretic definitions, in particular matchings in a bipartite graph.
Subsection \ref{sub4:sec2} relates polynomial matrices
and bipartite graphs; this subsection is relevant in the context of genericity
of parameters. Subsection \ref{sub5:sec2} contains a precise
definition of genericity and some simple examples. 

\subsection{Polynomial matrices}

Let $\R[s]$ be the commutative
ring of polynomials in the indeterminate
$s$ with coefficients from the field of real numbers $\R$. Let 
$\Rpvs$ be the ring of 
polynomial matrices with $p$ rows and $v$ columns  each of whose 
entry is a polynomial in $\Rs$.

A square polynomial matrix $M \in\Rvvs$ is said to be nonsingular if
$\det(M) \ne 0$. The roots of the polynomial $\det(M)$ are called the zeros
of $M$. Thus the zeros of a square polynomial matrix $M$ are the complex
numbers where $M$ loses its rank: we use this property to define
zeros of a possibly nonsquare polynomial matrix.
The zeros of $M(s)\in\Rpvs$ is defined to be the set of
those complex numbers $\lambda\in\C$
where the rank of the polynomial matrix `falls', more precisely:
\begin{equation} \label{zerosM}
{\zeros}(M):=\{\lambda \in \C ~|~ \rank(M(\lambda)) < \rank(M(s))  \}.
\end{equation}
The polynomial matrix $M\in\Rpvs$
is said to be full rank if $\rank(M)=\min(p,v)$. If $M$ is a full
rank polynomial matrix, $\zeros(M)$ can be found by computing the roots
of the greatest common divisor (gcd)
of the determinants of all the maximal minors of $M$. 
For a detailed exposition of these notions, we refer to \cite{Kai:1980}.

A polynomial matrix $U\in\Rvvs$ is called
unimodular if $\det(U)$ is a nonzero constant. These
are precisely the
square nonsingular polynomial matrices whose zero set is empty, or
equivalently, whose inverse is also a polynomial matrix.

\subsection{Behavioral approach}
A detailed exposition of the behavioural approach can be found in 
\cite{PolWil:1997}; we briefly cover 
the results that we need in this paper. A linear time invariant (LTI) 
dynamical system that is described by a system of ordinary differential
equations can be represented as
\begin{equation} \label{eq:R1-to-RN}
M_0w+M_1\der w+ \cdots M_N \deri{N}w=0 
\end{equation}
for constant matrices $M_i\in\R^{p\times v}$, with $M_N\ne 0$
and for $w$ a vector-valued, infinitely-often 
differentiable function $w:\R \rightarrow \R^{v}$. These $p$ equations
can be written in a shorthand notation by introducing the polynomial
matrix $M(s):=M_0+M_1s+\cdots +M_Ns^N$. Using the polynomial
matrix $M(s)\in \R^{p\times v}[s]$, the differential equations
in \eqref{eq:R1-to-RN} can
be written as $M(\der)w=0$. While the differential equations describing
a system are not unique, the set of trajectories that the system allows is
intrinsic to the system: we call the set of allowed trajectories as
{\em the behavior} of the system. More precisely, the allowed trajectories
are those that satisfy \eqref{eq:R1-to-RN}
\begin{equation} \label{eq:behavior}
\B=\{ w\in\CinfRRv ~|~ M(\der)w=0 \}.
\end{equation}
where $\CinfRRv$ denotes the space of infinitely often differentiable
functions from $\R$ to $\Rv$. The set $\B$ is called the behavior of the system.
In the context of \eqref{eq:behavior}, $M(\der)w=0$ is said
to be a {\em kernel representation} of the behavior.
A kernel representation is said to be 
minimal if the row dimension of $M$ is the minimum of all kernel representations
of $\B$: in this case $M$ has full row rank. This corresponds to the 
case that the equations describing the system are linearly independent over
$\R[s]$. 

A behavior
$\B$ is called controllable if for any $w_1,w_2\in\B$, there exist
$w_3\in\B$
and $T \in \R$ such that 
\[
\begin{array}{l}
w_3(t)=w_1(t)  \quad \mbox{ ~~  for $t<0$ ~~ and }\\
w_3(t)=w_2(t)  \quad \mbox{ ~~  for $t> T$. }
\end{array}
\]

A behavior $\B$
is called autonomous if one can conclude that $w_1=w_2$ whenever
$w_1,w_2\in\B$ satisfy $w_1(t)=w_2(t)$ for all $t\leqslant 0$.  
We state the required results from the behavioral literature
in the following proposition for easy
reference: see \cite{PolWil:1997, Wil:1991}.
\begin{prp} \label{prp:behavioral:results}
Consider $M\in\Rpvs$ and let behavior $\B$ be described
by the kernel representation $\Monw=0$. Then,
\begin{enumerate}
 \item $\B$ is autonomous if and only if $M$ has full column rank, i.e. 
    $\rank(M)=v$.
 \item The kernel representation $M$ is minimal if and only if $M$ has
  full row rank, i.e. $\rank(M)=p$.
 \item $\B$ is controllable if and only if $M(\lambda)$ has constant row rank
  for every complex number $\lambda \in \C$.
\end{enumerate}
\end{prp}
Thus, using the definition of $\zeros(M)$ as in 
\eqref{zerosM} above, a behavior described by $\Monw=0$ is controllable 
if and only if the zero set of $M$ is empty.  We use this characterization
of controllability and give equivalent graph theoretic conditions under
the assumption of genericity of parameters.

\subsection{Matchings in a bipartite graph} \label{sub3:sec2}

A graph $G=(V,E)$ in which $V$ can be partitioned into two non-empty 
sets $R$ and $C$ such that each edge in $E$ is 
between a vertex in $R$ and a vertex in
$C$ is called a bipartite graph. We use $G=(R,C;E)$ to indicate
these two vertex sets and the edge set.
For this paper, one vertex set $R$ denotes the rows and the other $C$ denotes
columns of the polynomial matrix $M$ describing the differential-algebraic
equations of \eqref{eq:behavior}.
A subgraph in which every vertex has degree at most one is
called a matching, i.e. each vertex has at most one edge
of this subgraph incident on it. The
number of edges in a matching $M$ is denoted by $|M|$.
For a bipartite graph $G=(R,C;E)$ with vertex sets $R$ and $C$, a matching
$m$ is said to be an {\em $R$-saturating matching} if $|m| = |R|$.
We define a $C$-saturating matching analogously. The special
case when $G$ satisfies $|R|=|C|$, an $R$-saturating matching
is also $C$-saturating matching, and vice-versa: we call such a matching
a {\em perfect matching}.  A detailed exposition on
these notions can be found in \cite{LovPlu:1986}.

A matching $M$ corresponds to a subset $r$ of $R$, $c$ of $C$ and
edges $e \subseteq E$.
We use $R(M), C(M)$ and $E(M)$ to denote subset $r$ of $R$,
$c \subseteq C$ and
the subset of edges that occur in the matching $M$ respectively.

As an example, consider the matrix $M\in\Rpvs$ with its nonzero
entries: $e_{ij}$, marked by its row and column indices $i$ and $j$ respectively:
\begin{equation} \label{eq:example}
M:=\left[ \begin{array}{ccc} 0 & e_{12} & e_{13} \\
e_{21} & e_{22} & e_{23}\end{array} \right].
\end{equation}

\comment{
\begin{figure}[!h]
\centerline{\resizebox{!}{45mm}{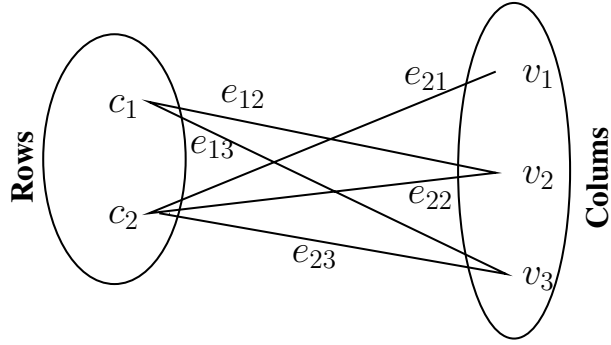}}
  \caption{Example of a bipartite graph}
  \label{bipart.expl.fig}
\end{figure}
}

In the associated bipartite graph (see Figure \ref{bipart.expl.fig}),
there are four $R$-saturating
matchings, say, $M_1$, $M_2$,
$M_3$ and $M_4$ corresponding to the edge pairs: $\{e_{12},e_{21} \}$, 
$\{e_{12},e_{23} \}$,
$\{e_{13},e_{22} \}$ and $\{e_{13},e_{21} \}$, respectively.

\subsection{Bipartite graph associated to polynomial matrices} 
\label{sub4:sec2}

We frequently need to associate an edge weighted bipartite graph to a 
given polynomial
matrix and vice-versa. For $M\in \Rpvs$, we associate an edge weighted 
bipartite
graph $G=(R,C;E)$ as follows. We set $R$ the rows of $M$ and $C$ as the
columns of $M$ and an edge between vertex $v_i$ of $R$ and $v_j$ of $C$
if the $(i,j)$-th entry in $M$ is nonzero. Further, the
degree of the polynomial in the $(i,j)$-the entry is assigned as
the weight of this edge. Thus $G$ has as many edges as the number of
nonzero entries in $M$ and the weights of the edges are non-negative integers.
Note that zero weight does {\em not} mean no edge; it means that entry in
the matrix $M$ is a nonzero constant, i.e. a polynomial of degree zero.
(The zero polynomial is usually said to be of degree $-\infty$; we do not
need this in our paper.) We say $G=(R,C;E)$ is the bipartite graph
associated to $M$. Conversely, given a bipartite graph $G=(R,C;E)$ with
the edges having non-negative integral weights, we can associate polynomial
matrices that satisfy the degree specifications in $E$ to the graph $G$.

We need some facts about the relation between the determinants of maximal minors
in a polynomial matrix $M$ and
$R$-saturating matchings of the bipartite graph $G$ associated
to $M$.  An elaborate treatment of these issues can be found in 
\cite{BabaiFrankl:1992}.
We first assume $M$ is square: suppose $M \in \Rpps$ is a polynomial matrix.
Let $m$ be a perfect matching in $G$.
Then $m$ corresponds to a nonzero term in the determinant expansion of $M$.
The determinant expansion of $M$ 
is the sum over all perfect matchings in $G$ (with suitable signs).
The matrix $M$  being nonsingular implies that there
exists at least one perfect matching in $G$. For the converse,
i.e., to conclude that $M$ is nonsingular if there are one or more
perfect matchings, one needs the important assumption of `genericity'.
While this is elaborated in the following subsection, we note here
that if there is only one perfect matching in $M$, then $M$ is nonsingular,
since no cancellations can occur. 
An upper triangular, lower triangular, or a diagonal square matrix has only
one contributing
perfect matching,
and hence the determinant comprises of just one
product and is nonzero (assuming the entries along the diagonal are nonzero).

We now deal with the case that $M\in\Rpvs$ is not square. Suppose
$p\leqslant v$ and there exists at least one $R$-saturating matching 
in the graph $G=(R,C;E)$ associated to $M$. (The case when
there does not exist even one $R$-saturating matching is dealt in Section
\ref{sec8} as `non-minimal descriptions'. If
$p\geqslant v$ then we will relate full column rank with $C$-saturating
matchings of $G$. This case is 
analogous to the $R$-saturating matching case by taking the transpose
of matrix $M$; this is not discussed further for this reason.)
Let $m$ be an $R$-saturating matching in $G$, and suppose $C(m)$ is
the subset of $C$ corresponding to this matching $m$. Then $m$ contributes
one term in the determinant expansion of the maximal minor corresponding
to $C(m)$ and $R$. In the situation that $M$ is not square, one
needs to deal with possibly different subsets $C(m_j)$ for different
$R$-saturating matchings $m_j$ in $G$. This is the crux of this paper and
is dealt with in the subsequent sections.

We return to the example of \eqref{eq:example} and Figure
\ref{bipart.expl.fig}. 
The maximal minor corresponding
to the first two columns has only one nonzero term: product of the entries
in the $R$-saturating matching $M_1$ defined after \eqref{eq:example}.
The maximal minor due to the
$2^{\rm nd}$ and $3^{\rm rd}$ columns has two nonzero terms,
namely products of the entries
in the matchings $M_2$ and $M_3$. Similarly, product of the entries
in the matching $M_4$ corresponds to the maximal minor: columns 1 and 3.

\subsection{Genericity of parameters} \label{sub5:sec2}

The notion of structural property makes the key assumption
of genericity of parameters. We use the following definition, as 
in \cite{Lin:1974} or \cite[page 132]{Mur:1987}.  
A set $S$ in $\Rn$ is said to be an algebraic variety
if $S$ is the solution set of an algebraic equation in $n$ variables.
The trivial equation
is the zero equation in the variables, in which case the algebraic variety
$S=\Rn$: we call this variety trivial.
We use the important fact that a nontrivial algebraic variety in $\Rn$
or $\Cn$ is a `thin' set, i.e. a set of measure zero. This is used
to define genericity of a property.

\begin{definition}\label{defn:genericity}
Consider property $P$ in terms of variables $a_1, \dots,a_{n}$ 
$\in\R$. Property
$P$ is said to be satisfied generically if the set of values $a_1,\dots$, $a_n$
that do {\em not} satisfy property $P$ form a nontrivial algebraic variety in $\Rn$.
\end{definition}

As a simple example, any two nonzero polynomials are generically
coprime. Let $M$ be a square matrix and suppose $G=(R,C;E)$ is its associated
bipartite graph. We saw that $M$ is nonsingular if $G$ has only one perfect
matching. If $G$ has two or more perfect matchings, then
the set of values of the parameters that cause cancellations of all 
terms, thus causing $M$ to be singular, is a so-called `thin' set, i.e.
these values form a non-trivial algebraic variety. More precisely,
these values form a set of measure zero in the space of all values that
can be attained by these parameters. Due to this reason,
we say $M$ is nonsingular generically if there exists
at least one perfect matching. A key assumption here is that the nonzero entries
in $M$ are chosen `independently' and hence they do not satisfy any 
nontrivial algebraic relation. In this context, the size of the maximal
matching on $G$ denotes the generic rank of $M$ and this rank is also
called the term-rank of $M$ (see \cite{Mur:1987}, for example).

\section{Polynomial matrices: generic properties}
\label{sec3}

In this section we state and prove
some generic properties of polynomial matrices; these are utilized
in the next section.  
The following lemma formulates necessary and sufficient conditions
for unimodularity of a polynomial matrix: it is one of the
main results of this paper.  The 
conditions for the matrix to be nonsingular is a standard
result from the literature, see \cite{Mur:1987}.

\begin{lemma} \label{thm:nonsingular:unimodular}
Consider the edge-weighted bipartite graph $G=(R,C;E)$ constructed 
from the polynomial matrix $M\in\Rvvs$.  
\begin{enumerate}
\item $M$ is generically nonsingular if and only if
the corresponding
bipartite graph has at least one perfect matching.
\item $M$ is generically unimodular if and only if 
     every perfect matching has only constant entry edges.
\end{enumerate}
\end{lemma}

\begin{proof1}
{\bf (1):} Since the determinant is a sum over all perfect 
matchings,
nonsingularity of $M$ implies existence of at least one perfect
matching in $G$.
Conversely, if one or more perfect matchings exist, then
we use
genericity to rule out cancellations and conclude that $M$ is nonsingular.
This proves Statement 1.

\noindent
{\bf (2):} If every perfect matching in $G$ has weight zero, then 
each product in the determinant expansion is a nonzero constant. Since cancellation
of these nonzero constants is ruled out upon addition, the determinant of
$M$ is also a nonzero constant, thus proving unimodularity. Conversely, if
there exists a perfect matching that has weight one or more, then
there is at least one product in the determinant expansion that has
degree one or more. Since cancellation is ruled out due to genericity,
the determinant of $M$ is of degree at least one, and hence
$M$ is not unimodular.  This proves Statement 2.  
\end{proof1}

We introduce the notion of a redundant edge: this plays a key role
in our results. Let $M\in\Rpvs$ with $p\leqslant v$ and 
consider the weighted bipartite graph $G=(R,C;E)$ constructed
from the polynomial matrix $M$.
An edge $e$ in $G$ is called {\em redundant}
if $e$ is not an element of any $R$-saturating matching of $G$. Thus the entry
corresponding to $e$ does not play a role in the determinant expansion
of any maximal minor of $M$;
this means $e$ does not affect the zero set of the polynomial matrix $M$.  
It turns out to help much in our results to remove the redundant edges
in a graph.
Consider the graph $G$, and remove all redundant
edges: we call
the resulting subgraph $G_{nr}$. It is the `maximal'
subgraph of $G$ with every edge non-redundant.
Clearly, $G$ has an $R$-saturating matching if and only if $G_{nr}$ has one.
Moreover, the removal of the non-redundant edges
results in a different polynomial matrix, say $M_{nr}$.
Due to the genericity assumption on $M$, and 
since the nonzero entries in $M_{nr}$ are same as those in $M$, 
we have the genericity property for $M_{nr}$ also.

Using $G_{nr}$,
the second statement in Lemma \ref{thm:nonsingular:unimodular}
can be restated
as follows.
The polynomial matrix $M$ is unimodular generically if and only if
all edges in $G_{nr}$ have weight zero.

\section{Structural controllability}
\label{sec4}

We now deal with the case that a polynomial matrix $M\in\Rpvs$ is not square.
Consider the bipartite graph $G=(R,C;E)$ associated to $M$.
The representation is generically minimal if and only if there exists
an $R$-saturating matching. Of course,
this requires $|R|\leqslant |C|$, i.e. $p \leqslant v$.

After the definition of genericity in Subsection \ref{sub5:sec2} above,
we saw that two nonzero polynomials $p$ and $q$ are generically coprime.
This is just another way of stating that $M=[p~~q]$ has
its zero set empty generically for nonzero polynomials $p$ and $q$.
However, if degree of $p=0$ (i.e. $p$ is a nonzero constant), then $M$
has its zero set empty even when $q$ is allowed to be zero. These
two cases are formulated in more generality in
Theorem \ref{no:subset:thm} below.

More generally, we ask when does a polynomial matrix $M\in\Rpvs$ have
an empty zero set. We answer this question in Theorem \ref{no:subset:thm} below
after we see two examples. Consider polynomial matrices $M_1$ and
$M_2 \in \Rpvs$
below, in which the $\ast$'s denote nonzero entries, and $p \leqslant v$.

\begin{equation} \label{eq:upper:triangular}
M_1= \left[  \begin{array}{cccccc} p_{11} & \ast & \cdots & \ast & \cdots & \ast \\
0 & \ast & \cdots & \ast & \cdots & \ast \\ \vdots & \vdots & \ddots & \vdots & \cdots & \vdots \\
0 & \ast & \cdots & \ast & \cdots & \ast \end{array} \right], \qquad 
M_2= \left[  \begin{array}{cc} N_{11} & N_{12} \\
0 & N_{22}\end{array} \right]
\end{equation}

Here, $p_{11}$ is an arbitrary nonzero polynomial and
$N_{11}$ is square and nonsingular. Obviously,
$\zeros(p_{11}) \subseteq \zeros(M_1)$ and
$\zeros(N_{11}) \subseteq \zeros(M_2)$. Thus the existence of a square
nonsingular
block (after permutation of rows and columns, if necessary) within
an upper triangular block matrix means additional conditions on such
a block matrix are required to hold for the zero set of the matrix 
$M_1$ or $M_2$
to be empty. This additional condition is 
that the zero set of this block also be empty, in other words, that
this square and nonsingular block be unimodular.
The following theorem states that non-existence of such a block submatrix is
also sufficient for the zero set of the polynomial matrix to be empty, under
the genericity assumption. A nonsingular submatrix that forms the left-upper triangular
block  (after permutation of rows/columns, if necessary) is related to
the property that all $R$-saturating matchings in $G_{nr}$
take the row-set of this submatrix
to the same column subset of $C$.

\begin{theorem1} \label{no:subset:thm}
Consider the bipartite graph $G=(R,C;E)$ associated to the polynomial
matrix $M\in \Rpvs$. Assume $p \leqslant v$, i.e. $|R| \leqslant |C|$. Suppose
all the redundant edges in $G$ are removed to obtain $G_{nr}$. 
A necessary and sufficient condition for the zero set of $M$ to be empty
generically is as follows.
\begin{itemize}
\item[P:] If there exist subsets $r \subseteq R$ and $c \subseteq C$ with
$|r|=|c|$ such that every $R$-saturating matching $M$ matches
$r$ to $c$, then all the edges in $G_{nr}$ incident on $r$ have weight zero.
\end{itemize} 
\end{theorem1}

We note that there are two cases within
Property $P$. Either there do not exist subsets $r$ and $c$ such
that every $R$-saturating matching matches $r$ to $c$, or there do exist
such subsets $r$ and $c$, in which case the edges are required to
satisfy a weight condition. The proof makes a distinction between these
two cases.  Another point to note is that removal of redundant edges
from $G$ corresponds to removal of corresponding entries from $M$ to obtain,
say, $M_{nr}$.  Since the entries that have been removed do not
affect any of the maximal minors of $M$, the zero sets of $M$ and $M_{nr}$
are equal. In the context of genericity of parameters in $M_{nr}$, as noted
above,
since the nonzero entries in $M_{nr}$ are the same as those in $M$, 
we have the genericity property for $M_{nr}$ also.

\begin{proof1}
{\bf (Necessity:)} Suppose there exist subsets $r$ and $c$ such that
$|r|=|c|$ and every $R$-saturating matching matches
$r$ to $c$, but there are some edges incident
on $r$ that have a nonzero weight. After a permutation
of elements of $R$ and $C$, the matrix $M$ is now in the form $M_2$ of
\eqref{eq:upper:triangular}, with the topmost $|r|$ rows 
corresponding to $r$ and leftmost $|r|$ columns corresponding to $c$.
Further, due to genericity, since one or more edges incident on $r$ have
weight at least one, the square block $N_{11}$ has determinant
a polynomial of degree at least one. This implies that at the roots
of the determinant, the matrix $N_{11}$ and hence $M_2$ loses rank,
thus showing that
the zero set of $M_2$ cannot be empty. This proves necessity of
the property $P$.

\noindent
{\bf (Sufficiency:)}
We now assume that the property $P$ is true, and show that the zero
set of $M$ is empty generically. We prove this by induction on
$|R|$. Let $|R|=1$. Since $P$ is true, either there exists a subset
$c\subseteq C $ with $|c|=1$ and every $R$ saturating matching matches
$R$ to $c$, or there doesn't exist such a subset $c$.
In the former case, property $P$ forces the entry $M_{c}$ to be a nonzero
constant, and hence the zero set of $M$ is empty. In the latter case, 
there exist at least two sets (in fact, singleton sets) $c_1$ and $c_2$ 
with $M_{c_1} \ne 0$ and 
$M_{c_2} \ne 0$. By genericity, these two polynomials have no common roots, and
hence the zero set of $M$ is empty.

We now assume the sufficiency of property $P$ for the zero set of $M$ to
be empty when the size of $R \leqslant k$.

Let $M_{nr}\in\R^{(k+1) \times \bullet }[s]$ have $k+1$ rows and assume
$M$ satisfies the property $P$. Again, we first consider the case when
there exist sets $r_1 \subset R$ and $c_1 \subset C$ such that
$|r_1|=|c_1|$ and every $R$-saturating matching matches $r_1$ to $c_1$.
The case where $r_1 = R$ is not covered by the induction hypothesis, but
follows from statement 2 of Lemma \ref{thm:nonsingular:unimodular}.
After a permutation of the rows and columns of $M_{nr}$, we have $M_{nr}$ as follows
\[
M_{nr}=\left[ \begin{array}{cc} M_{r_1~c_1} & 0 \\
0 & M_{r_2~c_2} \end{array} \right] 
\]
where $r_2:=R-r_1$ and $c_2:=C-c_1$. 
We now note that $ \zeros(M) = \zeros(M_{r_1~c_1})\cup \zeros(M_{r_2~c_2})$.
Further, every $R$-saturating matching in $M_{nr}$ is a union of row-saturating
matchings of $M_{r_1~c_1}$ and $M_{r_2~c_2}$. Hence, assumption of property $P$
for $M_{nr}$ implies this property for $M_{r_1~c_1}$ and $M_{r_2~c_2}$ also.
Since $M_{r_1~c_1}$ and $M_{r_2~c_2}$ have at most $k$ rows,
by the induction hypothesis, they have empty zero sets. Hence $M_{nr}$ also
has an empty zero set. This proves the sufficiency for the case when
there exist subsets $r$ and $c$ that are matched to each other by every
$R$ saturating matching of $M_{nr}$.

Consider the other case when there do not exist subsets $r$ and $c$ of
$R$ and $C$ respectively such that every $R$ saturating matching
matches $r$ and $c$. Consider the gcd of all $(k+1)\times(k+1)$ minors
of $M_{nr}$. Due to the absence of subsets $r$ and $c$, the various
minors have no common
factor arising from the determinant of a fixed
square nonsingular submatrix. Due to genericity, there is no other
reason that can cause the minors to have a common root, and hence
the gcd is equal to $1$. Thus the zero set of $M_{nr}$ is empty.
This proves the sufficiency of property $P$ for the case that
$M_{nr}$ has $k+1$ rows, and by induction this proves the sufficiency
part of Theorem \ref{no:subset:thm}.

\end{proof1}

The existence of subsets $r$ and $c$ satisfying above property would
mean that the determinant of this submatrix is a factor of every maximal
minor. The requirement that nonredundant
edges incident on $r$ have weight zero ensures that this submatrix
is unimodular, i.e. the determinant is a nonzero constant.

It appears from the above theorem that for the zero set of $R$ to be
generically empty, one requires to check the necessary and sufficient condition
for every subset $r$ of $R$, thus suggesting {\em exponential}
running time.
However, the absence of redundant edges makes it easy to formulate
the existence/non-existence of such a subset in easily verifiable conditions.  
The following theorem relates the condition to the size of a connected
component of $G$. This theorem, another of the main results of this paper,
is one that allows use of standard graph theoretic algorithms to
check structural controllability: the algorithm is analyzed in the 
following section.

\begin{theorem1} \label{thm2}
Let $M\in\Rpvs$ be a polynomial matrix of full row rank.
Consider the bipartite graph $G=(R,C;E)$
constructed from the rows and columns of $M$ and assume
there exists at least one $R$-saturating matching.
Suppose
all redundant edges in the bipartite graph $G$ are removed to
obtain $G_{nr}$.
Let $g_1$, $g_2, \dots$ $g_c$ be the connected components of $G_{nr}$.
Then $M(\lambda)$ has full row rank for every
complex number $\lambda \in \C$ generically if and only if
\begin{itemize}
\item for every component $g_i$ satisfying
 \begin{equation} \label{eq:vertex:sets:equal}
 |R(g_i)|=|C(g_i)|,
 \end{equation}
all edges in $g_i$ have weight zero.
\end{itemize}
\end{theorem1}

\begin{proof1}
The proof becomes simpler if we permute the rows and columns
of the polynomial matrix $M$ such that the matrix assumes a simpler
form.
We permute the rows and columns of $M_{nr}$ such that
each connected component of $G_{nr}$ correspond to consecutive
rows/columns. Thus $M_{nr}$ is now in the form:
\[
M_{nr}=\left[ \begin{array}{cccc}
M_1 & 0   & \cdots & 0 \\
0   & M_2 & \cdots & 0 \\
\vdots & \vdots & \ddots & \vdots\\
0 & 0 & \cdots & M_c \end{array} \right]
\]
with $M_i$ the submatrices of $M_{nr}$ corresponding to the connected
components $g_i$. Moreover, $M_i$ is square if and only if
$|R(g_i)|=|C(g_i)|$. Since there exists an $R$-saturating matching in $M$,
and hence in $M_{nr}$, each $g_i$ satisfies $|R(g_i)|\leqslant|C(g_i)|$, and
there exists at least one row-saturating matching for each
component $g_i$. Further,
\[
\zeros(M_{nr})=\underset{i=1,\dots,c}{\cup} \zeros(M_i).
\]
With this background, we proceed to the proof.

\noindent
{\bf Only if part: } Suppose $M$ has full row rank for every complex
number generically, then we show that for every component $g_i$ of
$G_{nr}$ satisfying  $|R(g_i)|=|C(g_i)|$, we have all edges in $g_i$ have
weight zero. Suppose $g_i$ is such that $|R(g_i)|=|C(g_i)|$,
and one or more edges in $g_i$ have a nonzero weight.
Since each edge is non-redundant, the determinant of $M_i$ is a non-constant.
This results in the roots of $\det (M_i)$ causing the zero set of $M$ to
be non-empty, thus proving the necessity of the condition by contradiction.

\noindent
{\bf If part: } For this part, we need to show that each of the $M_i$
is such that its zero set is empty. There are two cases. \\
{\bf Case 1:} $M_i$ is such that $|R(g_i)|=|C(g_i)|$, or\\ 
{\bf Case 2:} $M_i$ is such that $|R(g_i)|<|C(g_i)|$. 

In the first case, by assumption, all edges in $g_i$ have weight
zero, and hence determinant of $M_i$ is a nonzero constant. This proves
that the zero set is empty.

For the second case, since there are multiple row-saturating
matchings in $g_i$, there are at least two nonsingular maximal minors in $M_i$.
Further, there is no entry that is common to all the terms across all
maximal minors: this follows because $g_i$ is connected and every edge
is non-redundant. By genericity, the gcd of the two or more maximal minors
is one. This proves that the zero set of $M_i$ is empty.
\end{proof1}

\section{An algorithm and its efficiency} \label{sec5}

\newcommand{\nonnegZ}{{\mathbb{Z}_{\geqslant 0}}}
This section contains an analysis of the running time of various
algorithms needed for checking structrual controllability of
a dynamical system by using Theorem \ref{thm2}.  We first present
the algorithm as pseudocode below.  We assume that the input matrix
has been mapped into a weighted bipartite graph $G=(R,C;E)$ with 
$|R|=p, |C|=v$ and weights $w:E \mapsto \nonnegZ$.
\footnote{We use $\nonnegZ$ to denote the set of non-negative integers.}
Doing this takes $O(pv)$ units of time.\\[2mm]


\noindent
{\bf Input:} A weighted bipartite graph $G = (R,C;E)$ with 
$|R| = p, |C|=v$ and weights $w:E \mapsto \nonnegZ$. 

\noindent
{\bf Output:} ``Structurally controllable'' if the system is structurally controllable 
and ``Structurally uncontrollable'' otherwise.

\begin{algorithmic}[1]
  \STATE{$G_{nr}$ = Remove\_redundant\_edges($G$)}
  \STATE{Let $A_1,A_2,\ldots,A_t$ be the connected components of $G_{nr}$.}\\
  \COMMENT{Comment: Let $A_i$ be a graph with vertex set $V(A_i)$ and edge 
  		set $E(A_i)$.} \\
  \COMMENT{Comment: Let $R(A_i) = R \cap V(A_i)$ and $C(A_i) = C \cap V(A_i)$.}
  \IF{all components $A_i$ with $|R(A_i)| = |C(A_i)|$ have $w(e) = 0  \> $
          for every $ e \in E(A_i)$ }
  		\PRINT \texttt{``System structurally controllable''}
	\ELSE
  		\PRINT \texttt{``System structurally uncontrollable''}
	\ENDIF
\end{algorithmic}

We analyze the running time of each step of the above algorithm.

\noindent
{\bf Step 1: Removal of redundant edges:} 
Recall that an edge $e$ is called redundant if $e$ is not contained 
in any $R$-saturating matching of $G$.
One way to remove redundant edges is by first labelling all edges
as redundant or non-redundant.
In order to label edge $e = \{x,y\}$ as redundant or non-redundant, 
consider the
subgraph $G'$ obtained by removing $e$, the two vertices $x,y$,
and find the size of the maximum cardinality matching in $G'$.
Since $G'$ is also bipartite with $R(G') = R - 1$, it is clear that 
if the maximum cardinality matching has $|R| - 1$ edges,
then $e$ is a non-redundant edge, and if a maximum cardinality saturating
matching has strictly less than $|R|-1$ edges, then $e$ is a redundant 
edge in $G$.  Do this labelling for all edges $e \in E$.
Finding a maximum cardinality matching in bipartite graphs can be done 
in $O(E \sqrt{V})$ time by finding appropriate augmenting paths.
We do not present details of this algorithm as it is standard, instead, we
refer the
reader to the algorithm of Hopcroft and Karp \cite{HopKarp73} 
(see also \cite[Page 696]{CLRS}).  Since this is done for each
edge $e \in E$, the running time of our algorithm for classification of
all edges as redundant or non-redundant takes $O(E^2 \sqrt{V})$ time.

\noindent
{\bf Step 2: Decomposition of $G_{nr}$ into its connected components:}
Once all redundant edges in $G$ have been removed, the algorithm
for decomposing $G_{nr}$ into its connected components is again
standard and can be done in $O(|E_{nr}|~\mathrm{\log}^*~(|R|+|C|))$ time
\cite[Page 522]{CLRS}, where
$|E_{nr}|$ is the number of edges in the
subgraph $G_{nr}$ obtained from $G$ after removing all its redundant edges.


\noindent
{\bf Steps 3-7: Connected component checking:}
For each connected component $A_i$ satisfying
$|R(A_i)|=|C(A_i)|$, it takes $|E(A_i)|$ operations to check the weights
of all edges. In other words, in at most these many operations, one
can determine whether or not $A_i$ corresponds to a unimodular submatrix.

\begin{lemma} \label{lem:runtime}
 Given a $p \times v$ matrix $M$, let $E$ be the number of non-zero entries 
 of $M$.  Let $\B$ be the behaviour of $M$. 
Then, there exists an algorithm taking $O(E^2 \sqrt{p+v})$ time to check 
if $\B$ is structurally controllable.
\end{lemma}
\begin{proof1} 
Using the steps listed above and the running time involved for individual
steps, it is clear that the running time of the
algorithm is at most $O(E^2 \sqrt{p+v}) + O(E \log^*(p+v)) + O(E)$ which
is $O(E^2 \sqrt{p+v})$, thus completing the proof.
\end{proof1}

It is evident that one requires significantly lesser operations 
than mentioned in Step 1 for classification of edges into
redundant and non-redundant: in the process of marking an 
edge $e_1$, if $e_1$ is non-redundant because it could be completed
to an $R$-saturating matching in $G$, then one 
marks all other edges in that $R$-saturating matching also
in $G$ as non-redundant.
Hence the number of edges remaining to be marked
is fewer in the next sweep, if $e_1$ is non-redundant. On the other hand,
if an edge is redundant, then its immediate removal from $G$
will quicken the procedure for marking 
other edges. It appears that for dense matrices, there will
be significant improvement by such multiple-marking/intermediate-removal
procedure.  Thus the complexity could be significantly better than
$O(E^2 \sqrt{V})$. A precise count of the run-time complexity is
an interesting problem worth exploring.

\section{State space systems} \label{sec6}

The results in this paper simplify determination of structural controllability
for the situation $\der x =Ax + Bu$, with $A\in\Rnn$ and $B\in\Rnm$. 
In this section we use the main results of this paper to obtain a novel
method for checking structural controllability of the regular state space
system.

Construct the polynomial matrix
$M(s):=\begin{bmatrix}s I -A & B \end{bmatrix}$. For the dimensions assumed
on $A$ and $B$, the polynomial matrix $M(s)\in\R^{n \times (n+m)}[s]$.
Construct the bipartite graph $G$ corresponding
to $M$ and remove all redundant edges to obtain $G_{nr}$. 
 The following theorem states
 that structural controllability of $(A,B)$ is equivalent to each state
in $G_{nr}$ being connected to some input vertex.

\begin{theorem1} \label{thm:state:space}
Consider the bipartite graph $G$ constructed for 
$M(s)=\begin{bmatrix}s I -A & B \end{bmatrix} \in 
\R^{n \times (n+m)}[s]$. Obtain $G_{nr}$ by removal of all redundant edges.
The system $(A,B)$ is generically controllable if and only if $G_{nr}$
has the property that each state is connected to some input vertex.
\end{theorem1}

Before we proceed to the proof, it is noteworthy that the above theorem
is just application of generic rank check to the Popov-Belevitch-Hautus (PBH)
test
for controllability. While methods in structural controllability literature have
used the rank condition on $[B~AB~\cdots~A^{n-1}B]$, the PBH test has
not been explored with as much depth.  When applying the generic rank
condition to the PBH test,
the removal
of redundant edges turns out to result in checking connectivity of
state vertices in the {\em undirected} graph $G_{nr}$. The classical methods
treat first order systems using a {\em directed} graph approach, and this
is the result of the view that each variable in the systems is either an
input or output; consequently,
the input is to be utilized to control the output in a desired fashion.
The behavioral approach allows studying control
without having to classify variables into inputs/outputs. This 
results in an undirected graph. Further, higher order dynamical systems
are as easily dealt as first order systems in this approach.

The proof of the above result becomes easier after some notation that
is relevant to the state space case.  
We index the $R$-vertex set of the graph
by $\bullettop{x}_i$ for $i=1,\dots,n$, while the $C$-vertex set has
vertices corresponding to states and inputs: $C$ is indexed by $x_i$ for
$i=1,\dots,n$ and $u_j$ for $j=1,\dots,m$.  Further,
since we are dealing with a system of first order differential equations,
all the edges have weight either zero or one. We use a thick edge
for an edge of weight one, and dotted edge for edge of weight zero.
Further, the thick edges are precisely the `parallel edges': namely the ones
that connect $\bullettop{x}_i$ to $x_i$, while the dotted edges are the
non-parallel edges: those that connect $\bullettop{x}_i$ to $x_j$ for 
$i\ne j$ and also those that connect $\bullettop{x}_i$ to $u_k$. In other words,
the parallel edges correspond to the diagonal entries in $M(s):=[sI-A~~B]$,
i.e. the degree one entries in $M(s)$.
There are exactly $n$ parallel edges, and these form one $R$-saturating
matching in $G$, and hence all the parallel edges are also in $G_{nr}$.
It is sufficiently many
non-parallel edges connecting $u_k$ to the states
that help controllability of $(A,B)$: this is the intuitive idea of
the proof.

\begin{proof2} {\bf Theorem \ref{thm:state:space}: }
Due to all the parallel edges being non-redundant, and due to each 
$\bullettop{x}_i$ being connected to $x_i$, we first infer that
vertices $\bullettop{x}_i$ and $x_i$ lie in
the same connected component $g$. 
Hence, for each connected component $g$, the condition 
$|R(g)|=|C(g)|$ is equivalent to the absence of any input vertex $u$
in $C(g)$. 

\noindent
{\bf (Only if part:)} We assume that there exists a state $x_i$ such
that $x_i$ is
not connected to any input vertex in $G_{nr}$, and show that $(A,B)$ is not controllable.
Consider $g$, the connected component of $G_{nr}$ which contains $x_i$.
Due to the fact that each $\bullettop{x}_j$ of $R(G_{nr})$ is connected
to $x_j$ (of $C(G_{nr})$), the assumption on $x_i$
implies that there is no input vertex
in $C(g)$. This implies that for $g$, we have $|R(g)|=|C(g)|$. This means
every $R$-saturating matching in $G_{nr}$ matches $R(g)$ and $C(g)$.
Since the parallel edges in $g$ have weight one each, the condition in
Theorem \ref{thm2} is not satisfied. Thus $(A,B)$ is not controllable.

\noindent
{\bf (If part:)} We now show that if every state $x_i$ in $G_{nr}$ is
connected to some input vertex $u_j$, then every connected component
$g_k$ of $G_{nr}$ satisfies the condition $|R(g_i)| < |C(g_i)|$; from
Theorem \ref{thm2} above, it then follows that the system $(A,B)$ is
generically controllable.

Consider the connected components of $G_{nr}$. Due to the $n$ parallel edges
in $G$ that connect each state $x_j$ and $\bullettop{x}_j$ forming
an $R$-saturating matching, we noted above that each parallel edge
is non-redundant, and hence in $G_{nr}$. Further, this also causes
$|R(g)|\leqslant |C(g)|$ for each component $g$ of $G_{nr}$.
Note that $|C(g)|-|R(g)|$ is precisely the number of input vertices 
in $g$. Thus a state $x_i$ is connected to some input vertex
if and only if 
$|R(g)|< |C(g)|$ for the component $g$ that contains $x_i$.
Hence assuming that $G_{nr}$ is such that every state is connected to
some input vertex implies that 
$|R(g)|< |C(g)|$ for every connected component of $G_{nr}$. By Theorem
\ref{thm2} above, this implies that $(A,B)$ is structurally controllable.
\end{proof2}

{\boldmath
 \begin{figure}[!h]
 \centerline{\resizebox{!}{65mm}{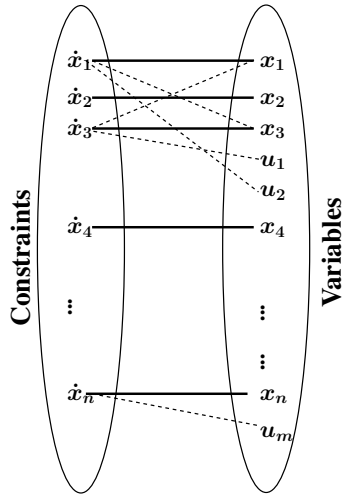}}
  \caption{State space case}
 \label{state.space.fig}
 \end{figure}
}

As shown in the above theorem, checking structural
controllability for a state space system reduces to checking
connectedness of every state vertex to some input vertex in the
bipartite graph $G_{nr}$.  Since the step classifying the edges in $G$ as
redundant or non-redundant (Step 1 of the algorithm in Section \ref{sec5})
is essential, and since this step was the most intensive of the steps
involved, the structural controllability question can
be answered in the state space case also in $O(E^2 \sqrt{m+n})$ time 
as shown in Lemma \ref{lem:runtime}.  
We use our method to determine structural controllability of some
state space models below.

\begin{example} \cite{Reinschke:1988}
Consider $A=\begin{bmatrix}0 & p & 0\\0 & q & 0\\ 0 & r & 0 \end{bmatrix}$
and $B=\begin{bmatrix}0 \\ s \\ 0\end{bmatrix}$,
for real numbers $p,q,r$ and $s$. This example has been shown to
be structurally
uncontrollable in \cite{Reinschke:1988}
using the method developed there.
By constructing the graph $G=(R,C;E)$ for $[sI-A~~B]$,
one can check that the edges corresponding to real numbers $p$ and
$r$ do not occur in any $R$-saturating matching, and hence $G_{nr}$
does not contain these edges. Consequently, $x_1$ and $x_3$ are not
connected to any input vertex.  
\end{example}

\begin{example}
Consider $A=\begin{bmatrix}0 & p & 0\\0 & q & 0\\ r & 0 & 0 \end{bmatrix}$
and $B=\begin{bmatrix}0 \\ s \\ 0\end{bmatrix}$,
for real numbers $p,q,r$ and $s$. This example has been obtained from
the previous one just by a change of location of entry $r$. The corresponding
graph $G_{nr}$ is now such that each state is connected to the
input vertex. 
\end{example}

\begin{example}
Consider $A=\begin{bmatrix}0 & r & 0\\0 & 0 & p\\ 0 & 0 & q \end{bmatrix}$
and $B=\begin{bmatrix}0 \\ 0 \\ s\end{bmatrix}$,
for real numbers $p,q,r$ and $s$. The significance of this example
is that it has been obtained from the previous one by just a permutation
of the states: this form is closer to the controller canonical form and
makes controllability of $(A,B)$ evident.
\end{example}

\section{Special models} \label{sec7}

In this section we study some special cases that are encountered often
in studying dynamical systems. We first deal with dynamical
systems whose models can be constructed using the fact that they
are made up of subsystems that have a particular property. We refer
to these as bond-graph-type systems and study them in Subsection \ref{sub7.1}.
The following subsection deals with analyzing the three most basic methods
of constructing a complex system: interconnection using series, parallel
and feedback connection. We prove a result about structural controllability
for these basic constructions. In Subsection \ref{sub7.3} we consider
the state space system in controller-canonical form and the Gilbert's 
canonical form and apply our main results to that case.

\subsection{Port-terminal interconnection based models} \label{sub7.1}

In this subsection, we consider an important class of mathematical models for
which we obtain finer bounds on the number of edges in the concerned bipartite
graph. These are systems comprised of subsystems that are interconnected
in a particular fashion. Many physical devices can be modeled as
having a certain `port-behavior' and such that when they are interconnected,
`energy' is exchanged across their ports. The interconnection of such
subsystems gives rise to larger systems that retain this property of
port-based energy exchange. These systems can also be analyzed using
bond-graph-based tools (\cite{MasSchBre:1991}).
 See also \cite{Wil:2010} for port-based
methods of studying this large class of physical systems.

The significance of such models for our paper is that
one can estimate the number
of edges in the polynomial matrix $M$ using the following heuristic arguments.
The matrix $M$ is composed of various differential/algebraic equations
corresponding to system laws. The system laws are typically of three
types, first: `device laws'. Each device
law has typically two or three variables. Secondly, interconnection laws
of the current/flow type (Kirchoff's current law, for example).
These are laws indicating that the net flow at
each junction is zero. Since each junction typically involves three or four
variables, the number of edges corresponding to these constraints also are
very few. The third important type of constraints are the `voltage drop' kind
of equations, this is akin to the Kirchoff's voltage laws. These equations arise from constraints that the net change of
`across variables' around each loop is zero; such constraints can involve
a fairly large number of variables, and these contribute to a large number
of edges. However, the number of such equations correspond to the number
of independent `loops', and hence these constraints are themselves typically
small in
number. Such arguments can be made in not just electrical networks, but
in any system that allows a bond-graph framework for modelling.

We assume $|E|\approx 3|R|$ to obtain heuristic estimates on the running 
time of the algorithm.  In this case, we get from Lemma \ref{lem:runtime}, that
the running time is $O(R^2 \sqrt{C})$. 

\subsection{Signal flow graphs}

In this subsection we consider the models for
dynamical systems constructed using signal-flow graphs.
We show that the three basic building blocks of complex
interconnections: the series, parallel and feedback connection
do not introduce any redundant edge in the resulting larger
mathematical model's bipartite graph $G$.
The absence of redundant edges results in significant improvement in
the runtime complexity of the algorithm proposed in Section \ref{sec5}: this
is elaborated before the beginning of Section \ref{sec8}.  
The larger question whether
arbitrary combination of the three building blocks still leads to:
1)~structural controllability, and 2)~no redundant edges, remains
an important open question.
%

The following theorem states that series, parallel and feedback
interconnection of two systems retains structural controllability, and
moreover, there are no redundant edges in the resulting bipartite graph.

\begin{theorem1} \label{thm:series:parallel}
Let $S_1$ and $S_2$ be two Single Input Single Output (SISO) systems.
Consider the system $S_3$
obtained by any one of the following interconnection procedures:\\
\hspace*{5mm}~~$\bullet$ series interconnection,\hspace*{1cm}
$\bullet$ parallel interconnection,\hspace*{1cm}
$\bullet$ feedback interconnection.\\
Then $S_3$ is structurally controllable and the bipartite graph constructed
for equations describing $S_3$ has no redundant edges.  
\end{theorem1} 

\begin{proof1}
Systems $S_1$ and $S_2$ are assumed to have transfer functions
$\frac{q_1(s)}{p_1(s)}$  and $\frac{q_2(s)}{p_2(s)}$  respectively.
We prove the result for 
the feedback interconnection, and give only the main features
for the other two interconnections.

Let $p_1(\der)y=q_1(\der)e$ and 
$p_2(\der)v=q_2(\der)y$ be the differential equations describing
$S_1$ and $S_2$. Feedback interconnection results in the additional
equation: $e=r-v$. Collating these three equations into a matrix, we get
$M(\der)w=0$ with 
\[
M_{\rm fdb}=\left[ \begin{array}{cccc}
  & 1 & 1 & -1 \\
p_1 & -q_1 & & \\
q_2 & & -p_2 & \end{array} \right] \mbox{ and } w_{\rm fdb}=
\left[ \begin{array}{c} y \\ e \\ v \\ r \end{array} \right] .
\]
The blank entries in the polynomial matrix $M_{\rm fdb}$ are all zero.
It is straightforward to see that each nonzero entry in $M_{\rm fdb}(s)$ 
occurs in some term of a suitable $3\times 3$ minor of $M_{\rm fdb}$. This means
that the bipartite graph constructed from $M_{\rm fdb}$
has no redundant edges, thus
proving the theorem for this interconnection configuration.  

For $S_1$ and $S_2$ connected in series and in parallel,
we can write the two systems of equations
in matrix form as follows:
\[
M_{\rm ser}=\left[ \begin{array}{ccc}
q_1 & -p_1 & \\
& q_2 & -p_2 \end{array} \right],~~ w_{\rm ser}=
\left[ \begin{array}{c} r \\ v \\ y \end{array} \right] 
\mbox{ and ~ }~~
M_{\rm par}=\left[ \begin{array}{cccc}
p_1 &  & -q_1 &  \\
& -p_2 & q_2 & \\ 
1 & 1 &   & -1 
\end{array} \right],~~ w_{\rm par}=
\left[ \begin{array}{c} u \\ v \\ r \\ y \end{array} \right] .
\]
Like the feedback interconnection case, non-redundancy of every edge
is verified by checking that each nonzero entry constitutes a term in
one or more maximal minors. This proves non-redundancy of every
edge. The structural controllability is verified by using
Theorem \ref{thm2}.  
\end{proof1}

\subsection{State space canonical forms} \label{sub7.3}

In this subsection we show that the
familiar controller canonical state space form
also has this property: there are no redundant edges in the bipartite
graph constructed from $M(s):=[sI-A~~B]$, when the $(A,B)$ pair is 
in the controller canonical form. We also show that
the Gilbert's canonical form also has no redundant edges
and displays structural controllability explicitly.

The next important situation when there are no redundant edges is
the familiar state space case when the pair $(A,B)$ is in the
controller canonical form.

\begin{theorem1}
Consider $A\in\Rnn$ and $B\in\Rn$ as in the equation below:
\[
A= \left[ \begin{array}{ccccc}
0 & 1 & 0 & \cdots & 0 \\
0 & 0 & 1 & \cdots & 0 \\
\vdots & \vdots & \vdots & \ddots & \vdots \\
0 & 0 & 0 & \cdots & 1 \\
-a_0 & -a_1 & -a_2 & \cdots & -a_{n-1} 
\end{array} \right] \mbox{ and  }
B= \left[ \begin{array}{c} 0 \\ 0 \\ \vdots \\ 0 \\ 1
\end{array} \right].
\]
Then the polynomial matrix $M(s):=[sI-A~~B]$ is structurally controllable
and has no redundant edges.  
\end{theorem1}

\begin{proof1}
The structural controllability is straightforward: from the matrix $M(s)$,
we see that there is one $n\times n$ minor whose determinant is equal to one;
namely the columns corresponding to the last $n$ columns. This proves
the full rank condition for every complex number $\lambda$ of $M(\lambda)$.

We now prove the absence of any redundant edge. Consider $M(s)=[sI-A~~B]$.
The diagonal entries in $M$ form a matching and hence are all non-redundant.
The non-redundance of the ones along the superdiagonal follows from the
previous paragraph. It remains to show that each of the $a_i$'s correspond
to non-redundant edges. This can be seen by expansion of the determinant
of $sI-A$ as follows. Since each of the $a_i$'s appear in the determinantal
expansion (in fact, they are the coefficients of the characteristic
polynomial), for each edge $a_i$, there is a maximal matching of size $n$
that contains this edge. This proves that these edges are also non-redundant.
This proves the theorem.
\end{proof1}

We now consider $A$ and $B$ in the Gilbert's canonical form (with 
$n=3$ for simplicity):
\[
A= \left[ \begin{array}{ccccc}
\lambda_1 & 1  &   \\
 & \lambda_1 &  \\
 & 0  & \lambda_3 
\end{array} \right] \mbox{ and  }
B= \left[ \begin{array}{c} 0 \\ b_2 \\ 0
\end{array} \right].
\]
It is a routine matter to verify
that each edge in the corresponding bipartite graph is non-redundant,
and that the state-space test for structural controllability
(Theorem \ref{thm:state:space}) allows conclusion of structural controllability
for this $(A,B)$ pair.

It is noteworthy that the controller
canonical form's realization (in terms of the classical series
interconnection of $n$ integrators; see \cite[page 39]{Kai:1980}, for example)
is just a combination of series and feedback configuration of several
SISO subsystems. An interesting open problem would be to prove that
there would be no redundant edges and one would have structural
controllability for arbitrary interconnection of SISO systems using
one or more of series, parallel and feedback interconnection
configurations.  In other words, it appears (and remains to be proved)
that the conventional signal flow graph satisfies the key properties of no
redundant edges and structural controllability.

The significance of the absence of redundant edges is that the
runtime complexity of the algorithm proposed in Section \ref{sec5}
can be significantly improved. Referring to that algorithm,
Step 1 would be inessential. Hence the running time is bounded by
$O(E\log^*(p+v)+E)$, which is $O(E\log^*(p+v))$, where recall that
$p$ and $v$ are the number of equations and variables, and $E$ is the
number of edges in $G$. Note that $G$ is same as $G_{nr}$.

\section{Non-minimal descriptions} \label{sec8}

One important case that we have not addressed so far was when
one or more equations describing the system are repeated or, more generally,
a linear combination of the other equations. This description
of the system is called non-minimal in behavioral literature. 
While manipulation of equations to obtain an equivalent minimal
set of equations is always possible when {\em exact} equations are specified,
this is not possible in the context of {\em structural} controllability
checks. This is because in checking structural controllability, we
only assume the structure of the system of equations is given, and further,
we make the key assumption that the parameters in the equations are
algebraically independent. This key assumption will fail to hold
if manipulation of equations is allowed to obtain a equivalent and
minimal description of the system. This section deals with
such non-minimal description of systems.

The first important point to note is that the rank of
a polynomial matrix is the size of a nonsingular minor of largest
size. If $M(\der)w=0$ is a kernel representation of a system, then
we are dealing with the case when $\rank(M) < $ row dimension of $M$.
Hence there does not exist an $R$-saturating matching in the bipartite
graph $G=(R,C;E)$ constructed from $M$. Suppose the size of the
matching with largest size is equal to $r_1$. Then controllability of
the system described by $M(\der)w=0$ is equivalent to coprimeness
of all the $r_1 \times r_1$ minors of $M$.  Since checking generic coprimeness
of all maximal minors is the main subject of this paper, the results
of this paper can easily be modified to handle the case of non-minimal
descriptions of linear dynamical systems.  The key modification
is that we now call an edge in $G$ redundant if it does not exist
in any matching of maximal size.

\begin{corollary}
Let $M(\der)w=0$ be a description, possibly non-minimal, of a dynamical
system. Construct the bipartite graph $G=(R,C;E)$ from $M$. Let $r_1$
be the size of the maximal matching in $G$. Construct $G_{nr}$ from
$G$  by removal of every redundant edge, i.e. an edge that doesn't
occur in any $r_1$ sized matching. Resolve $G_{nr}$ into its
connected components $g_i$. Then, the following are
equivalent.
\begin{itemize}
\item The system is structurally controllable,
\item For every component $g_i$ in $G_{nr}$
satisfying $ |R(g_i)|=|C(g_i)|$,
all edges in $g_i$ have weight zero.
\end{itemize}
\end{corollary}

The situation of non-minimal description of systems does
not happen in the state space case, and hence has not been addressed
in the literature. However, as noted in Section \ref{sec1},
this is relevant in the context of Smith
normal form of polynomial matrices: non-minimal description means
that the polynomial matrix does not have full rank.

\section{Conclusive remarks} \label{sec9}

We developed a method to check structural controllability of a system
of differential equations. As mentioned at the beginning of this
paper, this work is formulating conditions on a polynomial matrix $M$
under which
its invariant polynomials are generically one, i.e.
checking when the Smith normal form of $M$ has no
nonconstant polynomials along its diagonal. This related problem has
also been addressed for the case when the polynomial matrix does
not have full rank: the so-called non-minimal description of Section
\ref{sec9}. In addition to providing necessary and sufficient conditions
to check these properties, we also provided run-time
complexity of an algorithm to check this. For the more familiar state-space
desription of a dynamical system, this gives a novel method to check
controllability.  The central notion that we used at various edges
was that of a redundant edge: removal of redundant edges reveals
structural properties easily. It may be noted that redundancy here is
akin to the fact that the off-diagonal entries do not affect the
determinant of an upper-triangular (or lower-triangular) matrix.

In this paper, we addressed only the controllability aspect of
dynamical systems. The close relation between the methods to
check controllability and {\em observability} leave no reason to
address or pay any special attention to graph theoretic methods to
check structural {\em observability}. This is true for both
the state space and the behavioral description of dynamical
systems.

\vspace*{5mm}
\noindent
{\bf Acknowledgments:} We thank S.R~Khare, 
S.~Krishnan and D.~Chakraborty for useful discussions.

\bibliographystyle{IEEEtranS} 

\end{document}

%% file: example.pdftex_t
\begin{picture}(0,0)%
\includegraphics{example.pdf}%
\end{picture}%
\setlength{\unitlength}{3947sp}%
\begingroup\makeatletter\ifx\SetFigFont\undefined%
\gdef\SetFigFont#1#2#3#4#5{%
  \reset@font\fontsize{#1}{#2pt}%
  \fontfamily{#3}\fontseries{#4}\fontshape{#5}%
  \selectfont}%
\fi\endgroup%
\begin{picture}(4411,2516)(3555,-7304)
\put(7951,-6481){\rotatebox{90.0}{\makebox(0,0)[lb]{\smash{{\SetFigFont{14}{16.8}{\rmdefault}{\bfdefault}{\updefault}{\color[rgb]{0,0,0}Colums}%
}}}}}
\put(4276,-5611){\makebox(0,0)[lb]{\smash{{\SetFigFont{17}{20.4}{\rmdefault}{\bfdefault}{\updefault}{\color[rgb]{0,0,0}$c_1$}%
}}}}
\put(4276,-6436){\makebox(0,0)[lb]{\smash{{\SetFigFont{17}{20.4}{\rmdefault}{\bfdefault}{\updefault}{\color[rgb]{0,0,0}$c_2$}%
}}}}
\put(7321,-5386){\makebox(0,0)[lb]{\smash{{\SetFigFont{17}{20.4}{\rmdefault}{\bfdefault}{\updefault}{\color[rgb]{0,0,0}$v_1$}%
}}}}
\put(5101,-5536){\makebox(0,0)[lb]{\smash{{\SetFigFont{17}{20.4}{\rmdefault}{\bfdefault}{\updefault}{\color[rgb]{0,0,0}$e_{12}$}%
}}}}
\put(4876,-5911){\makebox(0,0)[lb]{\smash{{\SetFigFont{17}{20.4}{\rmdefault}{\bfdefault}{\updefault}{\color[rgb]{0,0,0}$e_{13}$}%
}}}}
\put(3721,-6286){\rotatebox{90.0}{\makebox(0,0)[lb]{\smash{{\SetFigFont{14}{16.8}{\rmdefault}{\bfdefault}{\updefault}{\color[rgb]{0,0,0}Rows}%
}}}}}
\put(6451,-5401){\makebox(0,0)[lb]{\smash{{\SetFigFont{17}{20.4}{\rmdefault}{\bfdefault}{\updefault}{\color[rgb]{0,0,0}$e_{21}$}%
}}}}
\put(5626,-6721){\makebox(0,0)[lb]{\smash{{\SetFigFont{17}{20.4}{\rmdefault}{\bfdefault}{\updefault}{\color[rgb]{0,0,0}$e_{23}$}%
}}}}
\put(6481,-6286){\makebox(0,0)[lb]{\smash{{\SetFigFont{17}{20.4}{\rmdefault}{\bfdefault}{\updefault}{\color[rgb]{0,0,0}$e_{22}$}%
}}}}
\put(7321,-6136){\makebox(0,0)[lb]{\smash{{\SetFigFont{17}{20.4}{\rmdefault}{\bfdefault}{\updefault}{\color[rgb]{0,0,0}$v_2$}%
}}}}
\put(7321,-6886){\makebox(0,0)[lb]{\smash{{\SetFigFont{17}{20.4}{\rmdefault}{\bfdefault}{\updefault}{\color[rgb]{0,0,0}$v_3$}%
}}}}
\end{picture}%

%% file: state.pdf_t
\begin{picture}(0,0)%
\includegraphics{state.pdf}%
\end{picture}%
\setlength{\unitlength}{3947sp}%
\begingroup\makeatletter\ifx\SetFigFontNFSS\undefined%
\gdef\SetFigFontNFSS#1#2#3#4#5{%
  \reset@font\fontsize{#1}{#2pt}%
  \fontfamily{#3}\fontseries{#4}\fontshape{#5}%
  \selectfont}%
\fi\endgroup%
\begin{picture}(4005,5960)(-14,-10431)
\put(676,-5236){\makebox(0,0)[lb]{\smash{{\SetFigFontNFSS{17}{20.4}{\rmdefault}{\bfdefault}{\updefault}{\color[rgb]{0,0,0}$\dot{x}_1$}%
}}}}
\put(676,-5686){\makebox(0,0)[lb]{\smash{{\SetFigFontNFSS{17}{20.4}{\rmdefault}{\bfdefault}{\updefault}{\color[rgb]{0,0,0}$\dot{x}_2$}%
}}}}
\put(676,-8236){\makebox(0,0)[lb]{\smash{{\SetFigFontNFSS{17}{20.4}{\rmdefault}{\bfdefault}{\updefault}{\color[rgb]{0,0,0}$\vdots$}%
}}}}
\put(3976,-8161){\rotatebox{90.0}{\makebox(0,0)[lb]{\smash{{\SetFigFontNFSS{20}{24.0}{\rmdefault}{\bfdefault}{\updefault}{\color[rgb]{0,0,0}Variables}%
}}}}}
\put(3001,-5236){\makebox(0,0)[lb]{\smash{{\SetFigFontNFSS{17}{20.4}{\rmdefault}{\bfdefault}{\updefault}{\color[rgb]{0,0,0}$x_1$}%
}}}}
\put(3001,-5686){\makebox(0,0)[lb]{\smash{{\SetFigFontNFSS{17}{20.4}{\rmdefault}{\bfdefault}{\updefault}{\color[rgb]{0,0,0}$x_2$}%
}}}}
\put(3001,-6061){\makebox(0,0)[lb]{\smash{{\SetFigFontNFSS{17}{20.4}{\rmdefault}{\bfdefault}{\updefault}{\color[rgb]{0,0,0}$x_3$}%
}}}}
\put(3001,-6436){\makebox(0,0)[lb]{\smash{{\SetFigFontNFSS{17}{20.4}{\rmdefault}{\bfdefault}{\updefault}{\color[rgb]{0,0,0}$u_1$}%
}}}}
\put(3001,-6811){\makebox(0,0)[lb]{\smash{{\SetFigFontNFSS{17}{20.4}{\rmdefault}{\bfdefault}{\updefault}{\color[rgb]{0,0,0}$u_2$}%
}}}}
\put(3001,-7261){\makebox(0,0)[lb]{\smash{{\SetFigFontNFSS{17}{20.4}{\rmdefault}{\bfdefault}{\updefault}{\color[rgb]{0,0,0}$x_4$}%
}}}}
\put(676,-7261){\makebox(0,0)[lb]{\smash{{\SetFigFontNFSS{17}{20.4}{\rmdefault}{\bfdefault}{\updefault}{\color[rgb]{0,0,0}$\dot{x}_4$}%
}}}}
\put(676,-6061){\makebox(0,0)[lb]{\smash{{\SetFigFontNFSS{17}{20.4}{\rmdefault}{\bfdefault}{\updefault}{\color[rgb]{0,0,0}$\dot{x}_3$}%
}}}}
\put(3001,-8311){\makebox(0,0)[lb]{\smash{{\SetFigFontNFSS{17}{20.4}{\rmdefault}{\bfdefault}{\updefault}{\color[rgb]{0,0,0}$\vdots$}%
}}}}
\put(3001,-8911){\makebox(0,0)[lb]{\smash{{\SetFigFontNFSS{17}{20.4}{\rmdefault}{\bfdefault}{\updefault}{\color[rgb]{0,0,0}$\vdots$}%
}}}}
\put(3001,-9286){\makebox(0,0)[lb]{\smash{{\SetFigFontNFSS{17}{20.4}{\rmdefault}{\bfdefault}{\updefault}{\color[rgb]{0,0,0}$x_n$}%
}}}}
\put(676,-9286){\makebox(0,0)[lb]{\smash{{\SetFigFontNFSS{17}{20.4}{\rmdefault}{\bfdefault}{\updefault}{\color[rgb]{0,0,0}$\dot{x}_n$}%
}}}}
\put(3001,-9736){\makebox(0,0)[lb]{\smash{{\SetFigFontNFSS{17}{20.4}{\rmdefault}{\bfdefault}{\updefault}{\color[rgb]{0,0,0}$u_m$}%
}}}}
\put(226,-8386){\rotatebox{90.0}{\makebox(0,0)[lb]{\smash{{\SetFigFontNFSS{20}{24.0}{\rmdefault}{\bfdefault}{\updefault}{\color[rgb]{0,0,0}Constraints}%
}}}}}
\end{picture}%